\documentclass[a4paper,12pt]{amsart}
\usepackage{amsfonts}
\usepackage{amssymb}
\usepackage{ifthen}
\usepackage{graphicx}
\nonstopmode \numberwithin{equation}{section}
\usepackage[usenames]{color}
\newtheorem{thm}{Theorem}[section]
\newtheorem{lem}{Lemma}[section]
\newtheorem{cor}[thm]{Corollary}
\newtheorem{prop}[thm]{Proposition}

\newtheorem{step}{Step}[section]

\theoremstyle{definition}
\newtheorem{mlem}{Main lemma}[section]
\newtheorem{assertion}{Assertion}[section]
\newtheorem{cl}{Claim}[section]
\newtheorem{ca}{Case}[section]
\newtheorem{sca}{Subcase}[section]
\newtheorem{scl}{Subclaim}[section]
\newtheorem{conj}[thm]{Conjecture}
\newtheorem{fact}{Fact}[section]
\newtheorem{defn}[thm]{Definition}
\newtheorem{op}[thm]{Open Problem}
\newtheorem{ques}[thm]{Question}
\newtheorem{rem}[thm]{Remark}
\newtheorem{exam}[thm]{Example}

\numberwithin{equation}{section}

\newcounter {own}
\def\theown {\thesection       .\arabic{own}}

\newenvironment{pf}[1][]{%
 \vskip 3mm
 \noindent
 \ifthenelse{\equal{#1}{}}%
  {{\slshape Proof. }}%
  {{\slshape #1.} }%
 }%
{\qed\bigskip}

\newcounter{alphabet}
\newcounter{tmp}
\newenvironment{Thm}[1][]{\refstepcounter{alphabet}%
\bigskip%
\noindent%
{\bf Theorem \Alph{alphabet}}%
\ifthenelse{\equal{#1}{}}{}{ (#1)}%
{\bf .} \itshape}{\vskip 8pt}

\makeatletter
\newcommand{\Ref}[1]{\@ifundefined{r@#1}{}{\setcounter{tmp}{\ref{#1}}\Alph{tmp}}}
\makeatother

\newenvironment{Lem}[1][]{\refstepcounter{alphabet}%
\bigskip%
\noindent%
{\bf Lemma \Alph{alphabet}}%
{\bf .} \itshape}{\vskip 8pt}

\newcounter{alphabet2}

\newcommand{\ID}{{\mathbb D}}

\newcommand{\diam}{{\operatorname{diam}}}

\newcommand{\area}{{\operatorname{Area}}}

\def\be{\begin{equation}}
\def\ee{\end{equation}}

\newcommand{\ben}{\begin{enumerate}}
\newcommand{\een}{\end{enumerate}}

\newcommand{\blem}{\begin{lem}}
\newcommand{\elem}{\end{lem}}
\newcommand{\bthm}{\begin{thm}}
\newcommand{\ethm}{\end{thm}}
\newcommand{\bcor}{\begin{cor}}
\newcommand{\ecor}{\end{cor}}
\newcommand{\beg}{\begin{exam}}
\newcommand{\eeg}{\end{exam}}
\newcommand{\begs}{\begin{examples}}
\newcommand{\eegs}{\end{examples}}
\newcommand{\bdefe}{\begin{defn}}
\newcommand{\edefe}{\end{defn}}
\newcommand{\bprob}{\begin{prob}}
\newcommand{\eprob}{\end{prob}}
\newcommand{\bques}{\begin{ques}}
\newcommand{\eques}{\end{ques}}
\newcommand{\bei}{\begin{itemize}}
\newcommand{\eei}{\end{itemize}}
\newcommand{\bcon}{\begin{conj}}
\newcommand{\econ}{\end{conj}}
\newcommand{\bop}{\begin{op}}
\newcommand{\eop}{\end{op}}

\newcommand{\bas}{\begin{assertion}}
\newcommand{\eas}{\end{assertion}}

\newcommand{\bfa}{\begin{fact}}
\newcommand{\efa}{\end{fact}}

\newcommand{\bca}{\begin{ca}}
\newcommand{\eca}{\end{ca}}

\newcommand{\bst}{\begin{step}}
\newcommand{\est}{\end{step}}

\newcommand{\bsca}{\begin{sca}}
\newcommand{\esca}{\end{sca}}

\newcommand{\bcl}{\begin{cl}}
\newcommand{\ecl}{\end{cl}}

\newcommand{\bmlem}{\begin{mlem}}
\newcommand{\emlem}{\end{mlem}}

\newcommand{\bscl}{\begin{scl}}
\newcommand{\escl}{\end{scl}}

\newcommand{\bcons}{\begin{conjs}}
\newcommand{\econs}{\end{conjs}}

\newcommand{\bprop}{\begin{prop}}
\newcommand{\eprop}{\end{prop}}

\newcommand{\br}{\begin{rem}}
\newcommand{\er}{\end{rem}}
\newcommand{\brs}{\begin{rems}}
\newcommand{\ers}{\end{rems}}
\newcommand{\bo}{\begin{obser}}
\newcommand{\eo}{\end{obser}}
\newcommand{\bos}{\begin{obsers}}
\newcommand{\eos}{\end{obsers}}
\newcommand{\bpf}{\begin{pf}}
\newcommand{\epf}{\end{pf}}
\newcommand{\ba}{\begin{array}}
\newcommand{\ea}{\end{array}}
\newcommand{\beq}{\begin{eqnarray}}
\newcommand{\beqq}{\begin{eqnarray*}}
\newcommand{\eeq}{\end{eqnarray}}
\newcommand{\eeqq}{\end{eqnarray*}}

\newcommand{\ds}{\displaystyle}

\newcounter{minutes}
\divide\time by 60
\newcounter{hours}
\multiply\time by 60 \addtocounter{minutes}{-\time}

\begin{document}

\bibliographystyle{amsplain}
\title []
{Lengths, area and modulus of continuity of some classes of complex-valued functions}

\def\thefootnote{}
\footnotetext{ \texttt{\tiny File:~\jobname .tex,
          printed: \number\day-\number\month-\number\year,
          \thehours.\ifnum\theminutes<10{0}\fi\theminutes}
} \makeatletter\def\thefootnote{\@arabic\c@footnote}\makeatother

\author{Shaolin Chen}
 \address{Sh. Chen, College of Mathematics and
Statistics, Hengyang Normal University, Hengyang, Hunan 421008,
People's Republic of China.} \email{mathechen@126.com}



\begin{abstract}
In this paper, we discuss the
modulus of continuity of solutions to Poisson's equation, and
 give bounds of length and area distortion for some classes of $K$-quasiconformal mappings satisfying  Poisson's equations.
 The obtained results are the extension of the corresponding classical results.
\end{abstract}

\maketitle \pagestyle{myheadings} \markboth{ Shaolin Chen}{Lengths, area and modulus of continuity of some classes of complex-valued functions}

\section{Preliminaries and  main results }\label{csw-sec1}

We use $\mathbb{C}$ to denote the complex plane. For
$a\in\mathbb{C}$ and  $r>0$,  let $\ID(a,r)=\{z:\, |z-a|<r\}$,
$\mathbb{D}_r=\mathbb{D}(0,r)$ and $\mathbb{D}=\mathbb{D}_1$, the
open unit disk in $\mathbb{C}$. Let $\mathbb{T}=\partial\mathbb{D}$
be the boundary of $\mathbb{D}$. Furthermore, we denote by
$\mathcal{C}^{m}(\Omega)$ the set of all complex-valued $m$-times
continuously differentiable functions from $\Omega$ into
$\mathbb{C}$, where $\Omega$ is a subset of $\mathbb{C}$ and
$m\in\{0,1,2,\ldots\}$. In particular,
$\mathcal{C}(\Omega):=\mathcal{C}^{0}(\Omega)$ denotes the set of
all continuous functions in $\Omega$.  Let $G$ be a domain of
$\mathbb{C}$, and let $\overline{G}$ be the closure of $G$.  We use $d_{G}(z)$ to denote the Euclidean distance
from $z$ to the boundary $\partial G$ of $G$. Especially, we always
use $d(z)$ to denote the Euclidean distance from $z$ to the boundary
of $\mathbb{D}$.

For a real $2\times2$ matrix $A$, we use the matrix norm
$$\|A\|=\sup\{|Az|:\,|z|=1\}$$ and the matrix function
$$\lambda(A)=\inf\{|Az|:\,|z|=1\}.$$

For $z=x+iy\in\mathbb{C}$, the
formal derivative of a complex-valued function $f=u+iv$ is given
by
$$D_{f}=\left(\begin{array}{cccc}
\ds u_{x}\;~~ u_{y}\\[2mm]
\ds v_{x}\;~~ v_{y}
\end{array}\right),
$$
so that
$$\|D_{f}\|=|f_{z}|+|f_{\overline{z}}| ~\mbox{ and }~ \lambda(D_{f})=\big| |f_{z}|-|f_{\overline{z}}|\big |,
$$
where $$f_{z}=\frac{1}{2}\big(
f_x-if_y\big)\;\;\mbox{and}\;\; f_{\overline{z}}=\frac{1}{2}\big(f_x+if_y\big).$$  We use
$$J_{f}:=\det D_{f} =|f_{z}|^{2}-|f_{\overline{z}}|^{2}
$$
to denote the {\it Jacobian} of $f$.

For $z, w\in\mathbb{D}$, let
\be\label{G-1}G(z,w)=\log\left|\frac{1-z\overline{w}}{z-w}\right|\ee
and

\be\label{P-1} P(z,e^{i\theta})=\frac{1-|z|^{2}}{|1-ze^{-i\theta}|^{2}}\ee
denote the {\it Green function} and  {\it (harmonic)
Poisson kernel}, respectively, where $\theta\in[0,2\pi]$.

Let $\psi:~\mathbb{T}\rightarrow\mathbb{C}$ be a bounded integrable
function  and let $g\in\mathcal{C}(\overline{\mathbb{D}})$. For
$z\in\mathbb{D}$, the solution to the {\it Poisson's equation}

$$\Delta f(z)=g(z)$$ satisfying the boundary condition
$f|_{\mathbb{T}}=\psi\in L^{1}(\mathbb{T})$ is given by

\be\label{eq-1.0} f(z)=P[\psi](z)-G[g](z),\ee
where
\be\label{eq-2.0}G[g](z)=\frac{1}{2\pi}\int_{\mathbb{D}}G(z,w)g(w)dA(w),~
~P[\psi](z)=\frac{1}{2\pi}\int_{0}^{2\pi}P(z,e^{it})\psi(e^{it})dt,\ee
and $dA(w)$ denotes the Lebesgue measure on $\mathbb{D}$. It is well
known that if $\psi$ and $g$ are continuous in $\mathbb{T}$ and in
$\overline{\mathbb{D}}$, respectively, then
$f=P[\psi]-G[g]$ has a continuous extension
$\tilde{f}$ to the boundary, and $\tilde{f}=\psi$ in $\mathbb{T}$
(see \cite[pp. 118-120]{Ho}  and \cite{K1,K2}).

A continuous increasing function $\omega:\, [0,\infty)\rightarrow
[0,\infty)$ with $\omega(0)=0$ is called a {\it majorant} if
$\omega(t)/t$ is non-increasing for $t>0$. Given a subset $\Omega$
of $\mathbb{C}$, a function $f:\, \Omega\rightarrow \mathbb{C}$ is
said to belong to the {\it Lipschitz space
$\mathcal{L}_{\omega}(\Omega)$} if there is a positive constant $C$
such that \be\label{eqd1} |f(z)-f(w)|\leq C\omega(|z-w|) ~\mbox{ for
all $z,\ w\in\Omega.$} \ee For $\delta_{0}>0$, let

\be\label{eqd2} \int_{0}^{\delta}\frac{\omega(t)}{t}\,dt\leq
C\cdot\omega(\delta),\ 0<\delta<\delta_{0}, \ee and
\be\label{eqd3}
\delta\int_{\delta}^{+\infty}\frac{\omega(t)}{t^{2}}\,dt\leq
C\cdot\omega(\delta),\ 0<\delta<\delta_{0}, \ee where $\omega$ is a
majorant and $C$ is a positive constant.

A majorant $\omega$ is said to be {\it regular} if it satisfies the
conditions (\ref{eqd2}) and (\ref{eqd3}) (see
\cite{Dy2,Dy1,P-1999}).

Let $G$ be a proper subdomain of $\mathbb{C}$. We say that a
function $f$ belongs to the {\it local Lipschitz space }
$\mbox{loc}\mathcal{L}_{\omega}(G)$ if (\ref{eqd1}) holds, with a fixed
positive constant $C$, whenever $z\in G$ and
$|z-w|<\frac{1}{2}d_{G}(z)$ (cf. \cite{GM,L}). Moreover, $G$ is said
to be a {\it $\mathcal{L}_{\omega}$-extension domain} if
$\mathcal{L}_{\omega}(G)=\mbox{loc}\mathcal{L}_{\omega}(G).$ The geometric
characterization of $\mathcal{L}_{\omega}$-extension domains was  given
by Gehring and Martio \cite{GM}. Then Lappalainen \cite{L}
generalized their characterization, and proved that $G$ is a
$\mathcal{L}_{\omega}$-extension domain if and only if each pair of points
$z,w\in G$ can be joined by a rectifiable curve $\gamma\subset G$
satisfying \be\label{eq1.0d}
\int_{\gamma}\frac{\omega(d_{G}(\zeta))}{d_{G}(\zeta)}\,ds(\zeta)
\leq C\omega(|z-w|) \ee with some fixed positive constant
$C=C(G,\omega)$, where $ds$ is the arc length measure on
$\gamma$.  Furthermore, Lappalainen \cite[Theorem 4.12]{L} showed
that $\mathcal{L}_{\omega}$-extension domains  exist only for majorants
$\omega$ satisfying  (\ref{eqd2}). 

The following result is the classical Hardy-Littlewood type Theorem for analytic functions  with respect to the majorant
$\omega(t)=\omega_{\alpha}(t)=t^{\alpha}~(0<\alpha\leq1)$ for $t\in[0,+\infty).$ In fact, the Hardy-Littlewood type Theorems
and the modulus of continuity of analytic functions are closely related.

\begin{Thm}{\rm (\cite[Theorem 5.1]{Du1})}\label{Th-Du-1}
Let $f$ be an analytic function in $\mathbb{D}$ and continuous in
$\overline{\mathbb{D}}$. Then
$$|f(e^{i\theta_{1}})-f(e^{i\theta_{2}})|\leq C\omega_{\alpha}(|\theta_{1}-\theta_{2}|)~\mbox{for
all}~0\leq\theta_{1},~\theta_{2}<2\pi$$ if and only if
$$|f'(z)|\leq C\frac{\omega_{\alpha}\big(d(z)\big)}{d(z)}~\mbox{for all}~z\in\mathbb{D},$$ where $C$ is a positive constant.
\end{Thm}

Krantz \cite{Kr} established the following Hardy-Littlewood type theorem
for real harmonic functions. 

\begin{Thm}{\rm (\cite[Theorem 15.8]{Kr})}\label{ThmA3}
Let $u$ be a real harmonic function in $\mathbb{D}$, and $\omega(t)=\omega_{\alpha}(t)=t^{\alpha}$ be a majorant for
$0<\alpha\leq1$. Then $u$ satisfies
$$|\nabla u(z)|\leq C\frac{\omega_{\alpha}\big(d(z)\big)}{d(z)}
~\mbox{ for all }~z\in\mathbb{D}
$$
if and only if
$$|u(z)-u(w)|\leq C\omega_{\alpha}(|z-w|)~\mbox{ for all }~z,w\in\mathbb{D},
$$ where $C$ is a positive constant.
\end{Thm}

Moduli of continuity of harmonic quasiregular mappings via
Hardy-Littlewood property is considered in \cite{AAM}. In \cite{MV},
the authors characterizes the moduli of continuity of a function $f$
by using the square of  distance function and module of $\Delta f$
(see the the class $OC^{2}(G)$ in \cite{MV}). In particular,
quasiregular versions of the well-known result due to Koebe,
\cite[Theorem 4.2]{M}, is established and, by using this result, an
extension of Dyakonov's theorem for quasiregular mappings in space
(without Dyakonov's hypothesis that it is a quasiregular local
homeomorphism), \cite[Theorem 4.3]{M}, is proved. The
charcterization of Lipschitz-type spaces for quasiregular mappings
by average Jacobian is also established in \cite[Theorem 4.3]{M}.

For a given $g\in\mathcal{C}(\overline{\Omega})$, let
$$\mathcal{F}_{g}(\Omega)=\{f\in\mathcal{C}(\overline{\Omega})\cap\mathcal{C}^{2}(\Omega):~\Delta
f(z)=g(z),~z\in\Omega\},$$ where $\Omega$ is a proper subdomain of
$\mathbb{C}$. Obviously, all analytic functions  and  harmonic
mappings defined in $\overline{\Omega}$ belong to
$\mathcal{F}_{0}(\Omega)$. We improve Theorems \Ref{Th-Du-1} and
\Ref{ThmA3} into the following form.

\begin{thm}\label{thm-1c} Suppose that $\omega$ is a majorant satisfying  \eqref{eqd2}, and $\Omega$ is a bounded
$\mathcal{L}_{\omega}$-extension domain. For a given $g\in\mathcal{C}(\overline{\Omega})$, let  $f\in\mathcal{F}_{g}(\Omega)$.
 Then $f\in\mathcal{L}_{\omega}(\Omega)$ if and only if there exists a constant $C>0$ such that,  for all $z\in\Omega$,

$$\|D_{f}(z)\|\leq C\frac{\omega\big(d_{\Omega}(z)\big)}{d_{\Omega}(z)}.
$$
\end{thm}

 A  mapping $f\in\mathcal{C}^{1}(\mathbb{D})$ is called a
{\it Bloch type mapping} if $f$ satisfies

$$\sup_{z\in\mathbb{D}}\left\{\|D_{f}(z)\|\omega\big((d(z))^{\alpha}\big)\right\}<+\infty,$$
where $\omega$ is a majorant and $\alpha>0$ is a constant. The set
of all Bloch type mappings, denoted by the symbol
$\mathcal{B}_{\omega}^{\alpha}$, forms a complex Banach space with
the norm $\|\cdot\|$ given by

$$\|f\|_{\mathcal{B}_{\omega}^{\alpha}}=|f(0)|+\sup_{z\in\mathbb{D}}\left\{\|D_{f}(z)\|\omega\big((d(z))^{\alpha}\big)\right\}.$$

In the following, by using  the weighted Lipschitz function,  Holland and  Walsh \cite{HW} gave an equivalent characterization of the analytic Bloch space.
For the related investigation of this topic for real
functions, we refer to \cite{Pav,Re}.

\begin{Thm}{\rm (\cite[Theorem 3]{HW})}\label{ThmA2} Let $f$ be
analytic in $\mathbb{D}$, and let $\omega$ be a majorant satisfying
$\omega(t)=t$ for $t\in[0,+\infty).$ Then
$f\in\mathcal{B}_{\omega}^{1}$ if and only if
$$\sup_{z,w\in\mathbb{D},z\neq w}\left\{
\frac{\sqrt{(1-|z|^{2})(1-|w|^{2})}|f(z)-f(w)|}{|z-w|}\right\}<\infty.
$$
\end{Thm}

In \cite{Dy1}, Dyakonov studied the relationship between the modulus of continuity and the bounded mean oscillation on analytic functions in $\mathbb{D}$, and obtained the following result.

\begin{Thm}{\rm (\cite[Theorem 1]{Dy1})}\label{ThmA}
Suppose that $f$ is an analytic function in $\mathbb{D}$ which is
continuous up to the boundary of $\mathbb{D}$. If $\omega$ and
$\omega^{2}$ are regular majorants,  then
$$f\in \mathcal{L}_{\omega}(\mathbb{D})\Longleftrightarrow P[|f|^{2}](z)-|f(z)|^{2}\leq
M\omega^{2}(d(z)).
$$

\end{Thm}

Analogy Theorems \Ref{ThmA2} and \Ref{ThmA}, we prove the following result.

\begin{thm}\label{thm-c3.0}

For a given $g\in\mathcal{C}(\overline{\mathbb{D}})$, let  $f\in\mathcal{F}_{g}(\mathbb{D})$.
Then, for $1\leq\alpha< 2$ and  a majorant $\omega$, the following statements are equivalent:
\begin{enumerate}
\item[{\rm (1)}] $f\in\mathcal{B}_{\omega}^{\alpha};$

\item[{\rm (2)}] There exists a constant $C>0$ such that for all  $r\in(0,d(z)]$,
$$\frac{1}{|\mathbb{D}(z,r)|}\int_{\mathbb{D}(z,r)}\left|f(\zeta)-
f(z)\right|dA(\zeta)\leq C\frac{r}{\omega(r^{\alpha})},
$$ where $|\mathbb{D}(z,r)|$ denotes  the area of $\mathbb{D}(z,r)$.
\end{enumerate}
\end{thm}

By \cite[Theorem 3]{CPR-2015} and Theorem \ref{thm-c3.0}, we obtain the following result which is a
generalization of Theorem \Ref{ThmA2}.

\begin{cor}
For a given $g\in\mathcal{C}(\overline{\mathbb{D}})$,
let  $f\in\mathcal{F}_{g}(\mathbb{D})$.
Then, for  $0\leq s<1$ and  $1\leq\alpha\leq s+1$, the following are equivalent:
\begin{enumerate}
\item[{\rm (1)}] $f\in\mathcal{B}_{\omega}^{\alpha};$

\item[{\rm (2)}] There exists a constant $C>0$ such that for all  $r\in(0,d(z)]$,
$$\frac{1}{|\mathbb{D}(z,r)|}\int_{\mathbb{D}(z,r)}\left|f(\zeta)-
f(z)\right|dA(\zeta)\leq\frac{Cr}{\omega(r^{\alpha})},
$$ where $|\mathbb{D}(z,r)|$ denotes  the area of $\mathbb{D}(z,r)$;

\item[{\rm (3)}] There exists a constant $C>0$ such that for all $z, w\in\mathbb{D}$ with $z\neq w$,
$$\frac{|f(z)-f(w)|}{|z-w|}\leq \frac{C}{\omega\big(d^{s}(z)d^{\alpha-s}(w)\big)}.
$$

\end{enumerate}
\end{cor}

For $r\in[0,1)$,   the {\it perimeter} of the curve
$C(r)=\big\{w=f(re^{i\theta}):\, \theta\in[0,2\pi]\big\}$, counting
multiplicity, is defined by

\be\label{eq-c-1}\ell_{f}(r)=\int_{0}^{2\pi}|df(re^{i\theta})|= r\int_{0}^{2\pi}\left|f_{z}(re^{i\theta})-e^{-2i\theta}
f_{\overline{z}}(re^{i\theta})\right|d\theta,
\ee
where  $f\in\mathcal{C}^{1}(\mathbb{D})$. In particular, let
$\ell_{f}(1)=\sup_{0<r<1}\ell_{f}(r)$ (cf. \cite{CLP-2017}).

A sense-preserving homeomorphic   $f$ from a domain $\Omega$ onto
$\Omega'$, contained in the {\it Sobolev class}
$W_{loc}^{1,2}(\Omega)$, is said to be a {\it $K$-quasiconformal
mapping} if, for $z\in\Omega$,
$$\|D_{f}(z)\|^{2}\leq K\big | \det D_{f}(z)\big |,~\mbox{i.e.,}~\|D_{f}(z)\|\leq K\lambda\big(D_{f}(z)\big),
$$
where $K\geq1$  (cf. \cite{K1,K2}). In the following, we will give bounds of length and area distortion for some classes of $K$-quasiconformal mappings satisfying  Poisson's equations.

\begin{thm}\label{thm-1}
For a given $g\in\mathcal{C}(\overline{\mathbb{D}})$,
let  $f\in\mathcal{F}_{g}(\mathbb{D})$. If $f=P[f]-G[g]$ is a $K$-quasiconformal
 mapping with $\ell_{f}(1)<+\infty$, then, for $n\geq1$,

\be\label{eq-ct-1}|a_{n}|+|b_{n}|\leq\frac{K\ell_{f}(1)}{2n\pi}+\frac{2}{3n}\|g\|_{\infty},\ee

\be\label{eq-ct-2}\sup_{z\in\mathbb{D}}\big\{\|D_{P[f]}(z)\|(1-|z|^{2})\big\}\leq\left(\frac{\ell_{f}^{2}(1)K}{4\pi^{2}}+\frac{4}{9}\|g\|_{\infty}^{2}
+\frac{\ell_{f}(1)K^{\frac{1}{2}}}{3\pi}\|g\|_{\infty}\right)^{\frac{1}{2}}\ee and
$f\in\mathcal{B}_{\omega}^{1},$
where $P[f](z)=\sum_{n=0}^{\infty}a_{n}z^{n}+\sum_{n=1}^{\infty}\overline{b}_{n}\overline{z}^{n}$ and $\omega(t)=t$.

In particular, if $K=1$ and $\|g\|_{\infty}=0$, then the estimates (\ref{eq-ct-1}) and (\ref{eq-ct-2}) are sharp, and the extreme function is $f(z)=z$ for $z\in\overline{\mathbb{D}}$.

\end{thm}

For $\theta\in[0,2\pi]$,   the {\it radial length} of the curve
$C_{\theta}(r)=\big\{w=f(\rho e^{i\theta}):\, 0\leq\rho\leq r<1\big\}$, counting
multiplicity, is defined by

\be\label{eq-cp-1}\ell_{f}^{\ast}(r,\theta)=\int_{0}^{r}|df(\rho e^{i\theta})|= \int_{0}^{r}\left|f_{z}(\rho e^{i\theta})+e^{-2i\theta}
f_{\overline{z}}(\rho e^{i\theta})\right|d\rho,
\ee
where  $f\in\mathcal{C}^{1}(\mathbb{D})$ (cf. \cite{CPR-Nolinear}). In particular, let
$$\ell_{f}^{\ast}(1,\theta)=\sup_{0\leq r<1}\ell_{f}^{\ast}(r,\theta).$$

\begin{thm}\label{thm-2}
For a given $g\in\mathcal{C}(\overline{\mathbb{D}})$,
let  $f\in\mathcal{F}_{g}(\mathbb{D})$. If $f=P[f]-G[g]$ is a $K$-quasiconformal
 mapping with $M=\sup_{\theta\in[0,2\pi]}\ell_{f}^{\ast}(\theta,1)<+\infty$,
then \be\label{eq-cs} |a_{n}|+|b_{n}|\leq
KM+\frac{2}{3}\|g\|_{\infty}~\mbox{ for $n\geq 1$,} \ee where
$P[f](z)=\sum_{n=0}^{\infty}a_{n}z^{n}+\sum_{n=1}^{\infty}\overline{b}_{n}\overline{z}^{n}$.
In particular, if $K=1$ and $\|g\|_{\infty}=0$, then the estimate
\eqref{eq-cs} is sharp and the extreme function is $f(z)=Mz$.
\end{thm}

The proofs of Theorems \ref{thm-1c}$\sim$\ref{thm-2} will be presented in
Section \ref{csw-sec2}.




\section{The proof of the  main results }\label{csw-sec2}

The following result easily follows from \cite[Lemma 2.7]{K2}.

\begin{Lem}\label{Lem-1}
 If  $g\in\mathcal{C}(\overline{\mathbb{D}})$, then, for $z\in\mathbb{D}$,
 $$\max\left\{\left|\frac{\partial}{\partial z}G[g](z)\right|,~\left|\frac{\partial}{\partial \overline{z}}G[g](z)\right|\right\}\leq\frac{1}{3}\|g\|_{\infty},$$
 where $G[g]$ is defined in {\rm (\ref{eq-2.0})}.
\end{Lem}

\subsection*{Proof of Theorem \ref{thm-1c}}  We first prove the necessity. Let $z\in\Omega$ and $r=d_{\Omega}(z)/2$. For $w\in\mathbb{D}(z,r)$, we have

\beqq
f(w)
=J_{1}(w)-J_{2}(w),
\eeqq
where $$J_{1}(w)=\frac{1}{2\pi}\int_{0}^{2\pi}P\left(\frac{w-z}{r},e^{i\theta}\right)f(z+re^{i\theta})d\theta$$ and $$J_{2}(w)=\frac{r^{2}}{2\pi}
\int_{\mathbb{D}}G\left(\frac{w-z}{r},\zeta\right)g(r\zeta+z)dA(\zeta),$$ where $G$ and $P$ are defined in {\rm (\ref{G-1})} and {\rm (\ref{P-1})}, respectively.
By elementary calculations,  we have

$$\frac{\partial }{\partial
w}G\left(\frac{w-z}{r},\zeta\right)=\frac{1}{2}\frac{r(|\zeta|^{2}-1)}{[r-(w-z)\overline{\zeta}](w-z-r\zeta)}$$

and

$$\frac{\partial }{\partial
\overline{w}}G\left(\frac{w-z}{r},\zeta\right)=\frac{1}{2}\frac{r(|\zeta|^{2}-1)}{[r-(\overline{w}-\overline{z})\zeta](\overline{w}-\overline{z}-r\overline{\zeta})},$$
which give that

\beq\label{rr-t1} \|D_{J_{2}}(w)\|&=&\left|\frac{r^{2}}{4\pi}
\int_{\mathbb{D}}\frac{\partial }{\partial
w}G\left(\frac{w-z}{r},\zeta\right)g(r\zeta+z)dA(\zeta)\right|\\
\nonumber &&+\left|\frac{r^{2}}{4\pi}
\int_{\mathbb{D}}\frac{\partial }{\partial
\overline{w}}G\left(\frac{w-z}{r},\zeta\right)g(r\zeta+z)dA(\zeta)\right|\\
\nonumber
&\leq&\frac{r^{2}\|g\|_{\infty}}{4\pi}\int_{\mathbb{D}}\frac{1-|\zeta|^{2}}{\left|\frac{w-z}{r}-\zeta\right|\left|1-\frac{(w-z)\overline{\zeta}}{r}\right|}
dA(\zeta)\\ \nonumber &&+
\frac{r^{2}\|g\|_{\infty}}{4\pi}\int_{\mathbb{D}}\frac{1-|\zeta|^{2}}{\left|\frac{\overline{w}-\overline{z}}{r}-\overline{\zeta}\right|
\left|1-\frac{(\overline{w}-\overline{z})\zeta}{r}\right|}
dA(\zeta). \eeq By (\ref{rr-t1}), Lemma \Ref{Lem-1} and by letting
$\xi=\frac{w-z}{r}$, we see that

\beq\label{rr-t2}
\|D_{J_{2}}(w)\|\leq\frac{r^{2}\|g\|_{\infty}}{2\pi}\int_{\mathbb{D}}\frac{1-|\zeta|^{2}}{\left|\xi-\zeta\right|
\left|1-\xi\overline{\zeta}\right|}
dA(\zeta)\leq\frac{2}{3}\|g\|_{\infty}r^{2}.
\eeq

The elementary computations lead to

$$\frac{\partial }{\partial
w}P\left(\frac{w-z}{r},e^{i\theta}\right)=\frac{-(\overline{w}-\overline{z})|w-z-re^{i\theta}|^{2}-(r^{2}-|w-z|^{2})
(\overline{w}-\overline{z}-re^{-i\theta})}{|w-z-re^{i\theta}|^{4}},
$$
and
$$\frac{\partial }{\partial
\overline{w}}P\left(\frac{w-z}{r},e^{i\theta}\right)=\frac{-(w-z)|w-z-re^{i\theta}|^{2}-(r^{2}-|w-z|^{2})
(w-z-re^{i\theta})}{|w-z-re^{i\theta}|^{4}}.
$$
Then, for $w\in\mathbb{D}(z,r/2)$,

\beq\label{rr-t3}
\left|\frac{\partial }{\partial
w}P\left(\frac{w-z}{r},e^{i\theta}\right)\right|&\leq&\frac{|w-z|}{|w-z-re^{i\theta}|^{2}}+\frac{r^{2}-|w-z|^{2}}{|w-z-re^{i\theta}|^{3}}\\ \nonumber
&\leq&\frac{\frac{r}{2}}{\frac{r^{2}}{4}}+\frac{r^{2}}{\frac{r^{3}}{8}}=\frac{10}{r}
\eeq
and \be\label{rr-t4}\left|\frac{\partial }{\partial
\overline{w}}P\left(\frac{w-z}{r},e^{i\theta}\right)\right|\leq\frac{10}{r}.\ee

It follows from (\ref{rr-t2}), (\ref{rr-t3}) and (\ref{rr-t4}) that, for $w\in\mathbb{D}(z,r/2)$,

\beq\label{rr-t5}
\|D_{f}(w)\|&\leq&\left|\frac{1}{2\pi}\int_{0}^{2\pi}\frac{\partial}{\partial w}P\left(\frac{w-z}{r},e^{i\theta}\right)(f(z+re^{i\theta})-f(z))d\theta\right|\\ \nonumber
&&+\left|\frac{1}{2\pi}\int_{0}^{2\pi}\frac{\partial}{\partial \overline{w}}P\left(\frac{w-z}{r},e^{i\theta}\right)(f(z+re^{i\theta})-f(z))d\theta\right|\\ \nonumber
&&+\|D_{J_{2}}(w)\|\\ \nonumber
&\leq&\frac{10}{r\pi}\int_{0}^{2\pi}|f(z+re^{i\theta})-f(z)|d\theta+\frac{2}{3}\|g\|_{\infty}r^{2}.
\eeq

Since $f\in\mathcal{L}_{\omega}(\Omega)$, we know that there is a positive constant $C_{1}$ such that

\be\label{rr-t6}|f(z+re^{i\theta})-f(z)|\leq C_{1}\omega(r).\ee
Since $\Omega$ is a bounded domain, we see that there is a positive constant $C_{2}$ such that
\be\label{rr-t7} \frac{\omega(r)}{r}\geq \frac{\omega(\diam(\Omega))}{\diam(\Omega)}\geq \frac{2}{3}\|g\|_{\infty}C_{2}.\ee

By (\ref{rr-t5}), (\ref{rr-t6}) and (\ref{rr-t7}), we conclude that there is a positive constant $C$ such that
$$\|D_{f}(w)\|\leq C\frac{\omega(r)}{r}.$$

Next, we show that the
sufficiency. Since $\Omega$ is a
$\mathcal{L}_{\omega}$-extension domain, we see that for any $z_{1},z_{2}\in\Omega$,
by using (\ref{eq1.0d}), there is a rectifiable curve
$\gamma\subset\Omega$ joining $z_{1}$ to $z_{2}$ such that
\begin{eqnarray*}
|f(z_{1})-f(z_{2})|&\leq&\int_{\gamma}\|D_{f}(\zeta)\|\,ds(\zeta)
\leq C\int_{\gamma}\frac{\omega\big(d_{\Omega}(\zeta)\big)}{d_{\Omega}(\zeta)}\,ds(\zeta)
\leq C\omega(|z_{1}-z_{2}|)
\end{eqnarray*}
for some constant $C>0$. The proof of this theorem is complete.
\qed

\begin{lem}\label{lem-1.2} For a given $g\in\mathcal{C}(\overline{\mathbb{D}})$, let $f\in\mathcal{F}_{g}(\mathbb{D})$.
 Then, for $a\in\mathbb{D}$,
there is a positive constant $C$ such that

$$\|D_{f}(a)\|\leq\frac{1}{\pi r}\int_{0}^{2\pi}|f(a+re^{i\theta})-f(a)|d\theta+\frac{2\|g\|_{\infty}}{3}r,$$
where $r\in(0,1-|a|)$.
\end{lem}

\bpf For $z\in\mathbb{D}_{r}$, let $$F(z)=f(z+a)-f(a).$$ Then,  $z\in\mathbb{D}_{r}$, $$\Delta F(z)=\Delta f(z+a)=g(z+a).$$
By (\ref{eq-1.0}), we have

$$F(z)=\frac{1}{2\pi}\int_{0}^{2\pi}\frac{r^{2}-|z|^{2}}{|z-re^{i\theta}|^{2}}F(re^{i\theta})d\theta-\frac{r^{2}}{2\pi}\int_{\mathbb{D}}
\log\left|\frac{r-z\overline{w}}{z-rw}\right|g(rw+a)dA(w)$$ for $z\in\mathbb{D}_{r}$.
By calculations, we have

\begin{eqnarray*}
F_{z}(z)&=&\frac{1}{2\pi}\int_{0}^{2\pi}\frac{-\overline{z}|z-re^{i\theta}|^{2}-(r^{2}-|z|^{2})(\overline{z}-re^{-i\theta})}{|z-re^{i\theta}|^{4}}F(re^{i\theta})d\theta\\
&&-\frac{r^{3}}{4\pi}\int_{\mathbb{D}}\frac{(|w|^{2}-1)}{(r-z\overline{w})(z-rw)}g(rw+a)dA(w)
\end{eqnarray*}
and

\begin{eqnarray*}
F_{\overline{z}}(z)&=&\frac{1}{2\pi}\int_{0}^{2\pi}\frac{-z|z-re^{i\theta}|^{2}-(r^{2}-|z|^{2})(z-re^{i\theta})}{|z-re^{i\theta}|^{4}}F(re^{i\theta})d\theta\\
&&-\frac{r^{3}}{4\pi}\int_{\mathbb{D}}\frac{(|w|^{2}-1)}{(r-w\overline{z})(\overline{z}-r\overline{w})}g(rw+a)dA(w),
\end{eqnarray*}
which 
yields that

\begin{eqnarray*}
\|D_{F}(0)\|&\leq&\frac{1}{r\pi}\int_{0}^{2\pi}|F(re^{i\theta})|d\theta+\frac{r}{2\pi}\int_{\mathbb{D}}\frac{(1-|w|^{2})}{|w|}|g(rw+a)|dA(w)\\
&\leq&\frac{1}{r\pi}\int_{0}^{2\pi}|F(re^{i\theta})|d\theta+\frac{2\|g\|_{\infty}}{3}r.
\end{eqnarray*}
The proof of this lemma is complete.
\epf

\subsection*{Proof of Theorem \ref{thm-c3.0}} We first prove $(1)\Rightarrow(2)$. By Lemma \ref{lem-1.2}, for
$\rho\in(0,d(z)]$,
$$\|D_{f}(z)\|\leq\frac{1}{\pi \rho}\int_{0}^{2\pi}\big|f(z+\rho e^{i\theta})-f(z)\big|\,d\theta+\frac{2\|g\|_{\infty}}{3}\rho,
$$
which gives
\beq\label{eq-rt-1}\int_{0}^{r}\rho^{2}\|D_{f}(z)\|d\rho&\leq&\frac{1}{\pi }\int_{0}^{r}\left(\rho\int_{0}^{2\pi}
|f(z+\rho e^{i\theta})-f(z)|d\theta\right)d\rho\\ \nonumber
&&+\frac{2\|g\|_{\infty}}{3}\int_{0}^{r}\rho^{3}d\rho,
\eeq
where $r=d(z).$ It follows from (\ref{eq-rt-1}) that
\begin{eqnarray*}
\|D_{f}(z)\|&\leq&\frac{3}{\pi
r^{3}}\int_{\mathbb{D}(z,r)}|f(z)-f(\zeta)|\,dA(\zeta)+\frac{\|g\|_{\infty}}{2}r\\
&=&\frac{3}{r |\mathbb{D}(z,r)|}\int_{\mathbb{D}(z,r)}|f(z)-f(\zeta)|\,dA(\zeta)+\frac{\|g\|_{\infty}}{2}r\\
&\leq&\frac{3C}{\omega(r^{\alpha})}+\frac{\|g\|_{\infty}}{2}r,
\end{eqnarray*} which gives that $f\in\mathcal{B}_{\omega}^{\alpha}.$

Now we prove  $(2)\Rightarrow(1)$. Since
$f\in\mathcal{B}_{\omega}^{\alpha},$
we see that there is a positive constant $C$ such that
\be\label{eq-11f}
\|D_{f}(z)\|\leq\frac{C}{\omega(d^{\alpha}(z))}.
\ee
For $z\in\mathbb{D}$ and $\zeta\in\mathbb{D}(z,r)$, we have

$$\omega\big(d^{\alpha}(z+t(\zeta-z))\big)\geq\omega\left(\big(d(z)-t|z-\zeta|\big)^{\alpha}\right),~t\in[0,1],$$
which, together with
 (\ref{eq-11f}), yields that
\beq\label{eq-12f}
|f(z)-f(\zeta)|&\leq&|z-\zeta|\int_{0}^{1}\|D_{f}(z+t(\zeta-z))\|\,dt\\ \nonumber
&\leq&C|z-\zeta|\int_{0}^{1}\frac{dt}{\omega\big(d^{\alpha}(z+t(\zeta-z))\big)}\\ \nonumber
&\leq&C|z-\zeta|\int_{0}^{1}\frac{dt}{\omega\left(\big(d(z)-t|z-\zeta|\big)^{\alpha}\right)}\\ \nonumber
&=&C\int_{0}^{|z-\zeta|}\frac{dt}{\omega\left(\big(d(z)-t\big)^{\alpha}\right)}.
\eeq
By (\ref{eq-12f}), we conclude that
\beq\label{eq-13f}
\nonumber\frac{1}{|\mathbb{D}(z,r)|}\int_{\mathbb{D}(z,r)}|f(z)-f(\zeta)|\,dA(\zeta)
&\leq&\frac{C}{|\mathbb{D}_{r}|}\int_{\mathbb{D}_{r}}\left(\int_{0}^{|\xi|}
\frac{dt}{\omega\left(\big(d(z)-t\big)^{\alpha}\right)}\right)dA(\xi)
\\
&=&\frac{2C}{r^{2}}\int_{0}^{r}\rho\left(\int_{0}^{\rho}
\frac{dt}{\omega\left(\big(d(z)-t\big)^{\alpha}\right)}\right)d\rho.
\eeq
By exchanging integral order, we obtain
\beq\label{eq-14f}
\int_{0}^{r}\rho\left(\int_{0}^{\rho}
\frac{dt}{\omega\left(\big(d(z)-t\big)^{\alpha}\right)}\right)d\rho&=&\int_{0}^{r}\left(\int_{t}^{r}\rho
d\rho\right)\frac{dt}{\omega\left(\big(r-t\big)^{\alpha}\right)}\\ \nonumber
&\leq&r\int_{0}^{r}\frac{\big(r-t\big)^{\alpha}}{\omega\left(\big(r-t\big)^{\alpha}\right)}
\big(r-t\big)^{1-\alpha}\,dt\\ \nonumber
&\leq&\frac{r^{\alpha+1}}{\omega(r^{\alpha})}\int_{0}^{r}\big(r-t\big)^{1-\alpha}\,dt\\ \nonumber
&=&\frac{1}{2-\alpha}\frac{r^{3}}{\omega(r^{\alpha})}.
\eeq
It follows from (\ref{eq-13f}) and (\ref{eq-14f}) that

$$\frac{1}{|\mathbb{D}(z,r)|}\int_{\mathbb{D}(z,r)}|f(z)-f(\zeta)|\,dA(\zeta)\leq\frac{2C}{2-\alpha}\frac{r}{\omega(r^{\alpha})}.$$
 The proof of this theorem is
complete. \qed

The following result is well-known (cf. \cite{Ca}).

\begin{Lem}\label{Lem-A}
Among all rectifiable Jordan curves of a given length, the circle
has the maximum interior area.
\end{Lem}

 \subsection*{Proof of Theorem \ref{thm-1}} We first prove (\ref{eq-ct-1}). Since $P[f]$ is harmonic in $\mathbb{D}$, we see that  $\partial P[f](z)/\partial z$ and $\partial P[f](z)/\partial \overline{z}$ are analytic and anti-analytic, respectively.
Hence, by Cauchy's integral formula, we have

$$na_{n}=\frac{1}{2\pi i}\int_{|z|=r}\frac{\frac{\partial P[f](z)}{\partial z}}{z^{n}}\,dz~\mbox{ and }~
nb_{n}=\frac{1}{2\pi i}\int_{|z|=r}\frac{\overline{\left(\frac{\partial P[f](z)}{\partial \overline{z}}\right)}}{z^{n}}\,dz,$$
which, together with $\|D_{P[f]}\|\leq\|D_{f}\|+\|D_{G[g]}\|$, implies that

\beq\label{eqt1.1}
n\big(|a_{n}|+|b_{n}|\big)&=&\frac{1}{2\pi}\left( \Bigg|\int_{|z|=r}\frac{\frac{\partial P[f](z)}{\partial z}}{z^{n}}\,dz\Bigg|+
\Bigg|\int_{|z|=r}\frac{\overline{\left(\frac{\partial P[f](z)}{\partial \overline{z}}\right)}}{z^{n}}\,dz\Bigg|\right)\\ \nonumber
&\leq&\frac{1}{2\pi r^{n}}\int_{0}^{2\pi}r\|D_{P[f]}(re^{i\theta})\|d\theta\\ \nonumber
&\leq&\frac{1}{2\pi r^{n}}\int_{0}^{2\pi}r\big(\|D_{f}(re^{i\theta})\|+\|D_{G[g]}(re^{i\theta})\|\big)d\theta,
\eeq where $r\in(0,1)$.

By (\ref{eq-c-1}), we have

\beq \label{eqt1.2}
\ell_{f}(1)&\geq&\ell_{f}(r)=r \int_{0}^{2\pi}\left|f_{z}(re^{i\theta})-e^{-2i\theta}
f_{\overline{z}}(re^{i\theta})\right|d\theta\\ \nonumber
&\geq& r\int_{0}^{2\pi}\big(|f_{z}(re^{i\theta})|-|f_{\overline{z}}(re^{i\theta})|\big)d\theta\\ \nonumber
&\geq&\frac{r}{K}\int_{0}^{2\pi}\|D_{f}(re^{i\theta})\|\,d\theta.
\eeq
It follows from (\ref{eqt1.1}), (\ref{eqt1.2}) and Lemma \Ref{Lem-1}  that

\begin{eqnarray*}
n\big(|a_{n}|+|b_{n}|\big)&\leq&\frac{K\ell_{f}(1)}{2\pi r^{n}}+\frac{1}{2\pi r^{n}}\int_{0}^{2\pi}r\|D_{G[g]}(re^{i\theta})\|d\theta\\
&\leq&\frac{1}{2\pi r^{n}}\left(K\ell_{f}(1)+\int_{0}^{2\pi}\|D_{G[g]}(re^{i\theta})\|d\theta\right)\\
&\leq&\frac{1}{2\pi r^{n}}\left(K\ell_{f}(1)+\frac{4\pi}{3}\|g\|_{\infty}\right),
\end{eqnarray*}
which gives that
$$|a_{n}|+|b_{n}|\leq\inf_{r\in(0,1)}\Bigg[\frac{1}{2n\pi r^{n}}\left(K\ell_{f}(1)+\frac{4\pi}{3}\|g\|_{\infty}\right)\Bigg]
=\frac{K\ell_{f}(1)}{2n\pi}+\frac{2}{3n}\|g\|_{\infty}.$$

Next we prove (\ref{eq-ct-2}). Let $\area(f(\mathbb{D}_{r}))$ denote the area of $f(\mathbb{D}_{r}),$ where $r\in(0,1)$. Then

\be\label{eq-t4}
\area(f(\mathbb{D}_{r}))=\int_{\mathbb{D}_{r}}J_{f}(z)\,dA(z)\geq\frac{1}{K}\int_{\mathbb{D}_{r}}\|D_{f}(z)\|^{2}\,dA(z).
\ee
 For $\theta\in[0,2\pi]$ and $z\in\mathbb{D}$, let
$$H_{\theta}(z)=\frac{\partial P[f](z)}{\partial z}+e^{i\theta}\overline{\left(\frac{\partial P[f](z)}{\partial \overline{z}}\right)}.$$
Then, by the subharmonicity of $|H_{\theta}|^{2}$, we obtain

\beq\label{eq-t5}
|H_{\theta}(z)|^{2}&\leq& \frac{1 }{\pi(1-|z|^{2})^{2}}\int_{0}^{1-|z|^{2}}\rho\int_{0}^{2\pi}|H_{\theta}(z+\rho
e^{i\gamma})|^{2} \,d\gamma \,d\rho\\ \nonumber
&\leq&\frac{1}{\pi(1-|z|^{2})^{2}}\int_{\mathbb{D}_{1-|z|^{2}}}\|D_{P[f]}(z+\zeta)\|^{2}\,dA(\zeta)\\ \nonumber
&\leq&\frac{I}{\pi(1-|z|^{2})^{2}},
\eeq where $$I=\int_{\mathbb{D}}\big(\|D_{G[g]}(\xi)\|+\|D_{f}(\xi)\|\big)^{2}\,dA(\xi).$$

By (\ref{eq-t4}), Lemma \Ref{Lem-1} and Cauchy-Schwarz's inequality, we get

\beq\label{eq-t6}
I&=&\int_{\mathbb{D}}\|D_{f}(\xi)\|^{2}\,dA(\xi)+\int_{\mathbb{D}}\|D_{G[g]}(\xi)\|^{2}\,dA(\xi)\\ \nonumber
&&+2\int_{\mathbb{D}}\|D_{f}(\xi)\|\|D_{G[g]}(\xi)\|\,dA(\xi)\\ \nonumber &\leq&\int_{\mathbb{D}}\|D_{f}(\xi)\|^{2}\,dA(\xi)+\frac{4\pi}{9}\|g\|_{\infty}^{2}+\frac{2}{3}\|g\|_{\infty}
\int_{\mathbb{D}}\|D_{f}(\xi)\|\,dA(\xi)\\ \nonumber
&\leq&K\area(f(\mathbb{D}))+\frac{4\pi}{9}\|g\|_{\infty}^{2}\\ \nonumber &&+\frac{2}{3}\|g\|_{\infty}\left(\int_{\mathbb{D}}\|D_{f}(\xi)\|^{2}\,dA(\xi)\right)^{\frac{1}{2}}
\left(\int_{\mathbb{D}}dA(\xi)\right)^{\frac{1}{2}}\\ \nonumber
&\leq&K\area(f(\mathbb{D}))+\frac{4\pi}{9}\|g\|_{\infty}^{2}+\frac{2\pi^{\frac{1}{2}}}{3}\|g\|_{\infty}\big(K\area(f(\mathbb{D}))\big)^{\frac{1}{2}}.
\eeq

Applying Lemma \Ref{Lem-A}, we have

$$\area(f(\mathbb{D}))\leq\pi\left(\frac{\ell_{f}(1)}{2\pi}\right)^{2}=\frac{\ell_{f}^{2}(1)}{4\pi},$$ which, together with (\ref{eq-t6}), yields that

\be\label{eq-t7}I\leq \frac{\ell_{f}^{2}(1)K}{4\pi}+\frac{4\pi}{9}\|g\|_{\infty}^{2} +\frac{\ell_{f}(1)K^{\frac{1}{2}}}{3}\|g\|_{\infty}. \ee

By (\ref{eq-t5}) and (\ref{eq-t7}), we conclude that

\be\label{eq-t17}\|D_{P[f]}(z)\|=\max_{\theta\in[0,2\pi]}|H_{\theta}(z)|\leq\frac{\left(\frac{\ell_{f}^{2}(1)K}{4\pi^{2}}+\frac{4}{9}\|g\|_{\infty}^{2}
+\frac{\ell_{f}(1)K^{\frac{1}{2}}}{3\pi}\|g\|_{\infty}\right)^{\frac{1}{2}}}{1-|z|^{2}}.\ee

At last, $f\in\mathcal{B}_{\omega}^{1}$ follows from (\ref{eq-t17}) and Lemma \Ref{Lem-1},
where $\omega(t)=t$.
The proof of this theorem is complete.
\qed

The following result is considered to be  a Schwarz-type lemma of subharmonic functions.

\begin{Thm} {\rm  (\cite[Theorem 2]{B})}\label{Thm-cs}
Let $\phi$ be subharmonic in $\mathbb{D}$. If, for all $r\in[0,1)$,
$$A(r)=\sup_{\theta\in[0,2\pi]}\int_{0}^{r}\phi(\rho e^{i\theta})\,d\rho\leq1,
$$
then $A(r)\leq r$.
\end{Thm}

\subsection*{Proof of Theorem \ref{thm-2}}

By Cauchy's integral formula, for $\rho \in (0,1)$ and $n\geq 1$, we get
$$na_{n}=\frac{1}{2\pi i}\int_{|z|=\rho}\frac{\frac{\partial P[f](z)}{\partial z}}{z^{n}}\,dz~\mbox{ and }~
nb_{n}=\frac{1}{2\pi
i}\int_{|z|=\rho}\frac{\overline{\left(\frac{\partial
P[f](z)}{\partial \overline{z}}\right)}}{z^{n}}\,dz,$$ which implies
that
\beq\label{eqc-30}
n(|a_{n}|+|b_{n}|)&=&\frac{1}{2\pi}\left|\int_{|z|=\rho}\frac{\frac{\partial
P[f](z)}{\partial z}}{z^{n}}\,dz\right|
+\frac{1}{2\pi}\left|\int_{|z|=\rho}\frac{\overline{\left(\frac{\partial
P[f](z)}{\partial \overline{z}}\right)}}{z^{n}}\,dz\right|\\
\nonumber &\leq&\frac{1}{2\pi
\rho^{n-1}}\int_{0}^{2\pi}\|D_{P[f]}(\rho e^{i\theta})\|\,d\theta.
\eeq
By calculations, for $\theta\in[0,2\pi]$, we obtain
\begin{eqnarray*}
\ell_{f}^{\ast}(\theta,r)
&=&\int_{0}^{r}|f_{z}(\rho e^{i\theta})+e^{-2i\theta}f_{\overline{z}}(\rho e^{i\theta})|\,d\rho\\
&\geq&\int_{0}^{r}\lambda(D_f)(\rho e^{i\theta})\,d\rho\\
&\geq&\frac{1}{K}\int_{0}^{r}\|D_{f}(\rho e^{i\theta})\|\,d\rho,
\end{eqnarray*}
which gives
\be\label{eqc-31}
\int_{0}^{r}\|D_{f}(\rho e^{i\theta})\|\,d\rho\leq K\ell_{f}^{\ast}(\theta,r)\leq KM.
\ee

It follows from   (\ref{eqc-31}) and Lemma \Ref{Lem-1} that

 \beq
 \int_{0}^{r}\|D_{P[f]}(\rho e^{i\theta})\| \, d\rho&\leq&\int_{0}^{r}\|D_{f}(\rho e^{i\theta})\|\,d\rho+\int_{0}^{r}\|D_{G[g]}(\rho e^{i\theta})\|\,d\rho\\ \nonumber
 &\leq&KM+\frac{2}{3}\|g\|_{\infty}r.
 \eeq

By (\ref{eqc-31}), the subharmonicity of $D_{P[f]}(\rho e^{i\theta})$ and Theorem \Ref{Thm-cs}, we have
\be\label{eqc-32}
\int_{0}^{r}\|D_{P[f]}(\rho e^{i\theta})\| \, d\rho\leq \left(KM+\frac{2}{3}\|g\|_{\infty}\right)r.
\ee
By (\ref{eqc-30}) and (\ref{eqc-32}), we get
\begin{eqnarray*}
2\pi n(|a_{n}|+|b_{n}|)\int_{0}^{r}\rho^{n-1}\, d\rho
&=&\int_{0}^{r}\left(\int_{0}^{2\pi}\|D_{P[f]}(\rho e^{i\theta})\|d\theta\right)\,d\rho\\
&=&\int_{0}^{2\pi}\left(\int_{0}^{r}\|D_{P[f]}(\rho e^{i\theta})\|d\theta\right) \, d\rho\\
&\leq&2\pi\left(KM+\frac{2}{3}\|g\|_{\infty}\right)r,
\end{eqnarray*}
which yields that
$$|a_{n}|+|b_{n}|\leq\inf_{r\in(0,1)}\left(\frac{KM+\frac{2}{3}\|g\|_{\infty}}{r^{n-1}}\right)=KM+\frac{2}{3}\|g\|_{\infty} ~\mbox{ for $n\geq  1$}.
$$ The proof of this theorem is complete.
\qed

\bigskip

{\bf Acknowledgements:}    This research was partly supported by the
Science and Technology Plan Project of Hengyang City (No.
2018KJ125), the National Natural Science Foundation of China (No.
11571216), the Science and Technology Plan Project of Hunan Province
(No. 2016TP1020), the Science and Technology Plan Project of
Hengyang City (No. 2017KJ183), and the Application-Oriented
Characterized Disciplines, Double First-Class University Project of
Hunan Province (Xiangjiaotong [2018]469).

\normalsize

\end{document}